\documentclass[a4paper,12pt]{amsart}
\usepackage{tikz-cd}
\usepackage{amssymb}
\usepackage{amsthm}
\usepackage{amsmath}
\usepackage{mathrsfs}

\newtheorem{theorem}{Theorem}[section]
\newtheorem{lemma}{Lemma}[section]
\newtheorem*{claim-1}{Claim 1}
\newtheorem*{claim-2}{Claim 2}
\newtheorem{proposition}{Proposition}[section]
\newtheorem{corollary}{Corollary}[section]
\newtheorem{theoremA}{Theorem}

\newtheorem{statementA}{Statement}

\theoremstyle{definition} 
\newtheorem{definition}{Definition}[section]
\newtheorem*{definition*}{Definition}

\newtheorem{problem}{Problem}[section]
\newtheorem{problemA}{Problem}

\theoremstyle{remark} 
\newtheorem{remark}{Remark}[section]

\def\bigdiag{\operatornamewithlimits{%
  \mathchoice{\vcenter{\hbox{\Large$\Delta$}}}
             {\vcenter{\hbox{\large$\Delta$}}}
             {\Delta}
             {\Delta}}}
\def\diag{\operatornamewithlimits{%
  \mathchoice{\vcenter{\hbox{$\Delta$}}}
             {\vcenter{\hbox{$\Delta$}}}
             {\Delta}
             {\Delta}}}

\begin{document}

{
\renewcommand*{\thefootnote}{$\star$}

\title{$\mathbb R^{\omega_1}$-Factorizable Spaces and Groups}
\footnotetext[0]{The work of O.~Sipacheva was financially supported by the Russian Science Foundation, 
grant~22-11-00075-P.} }

\author{Anton Lipin}
\email[A. Lipin]{tony.lipin@yandex.ru}
\address{Krasovskii Institute of Mathematics and Mechanics, ul.~S.~Kovalevskoi 16,  Yekaterinburg,
            620990 Russia \\
Ural Federal University, ul.~Mira 19, Yekaterinburg, 620002 Russia}

\author{Evgenii Reznichenko}
\email[E. Reznichenko]{erezn@inbox.ru} 
\address{Department of General Topology and Geometry, Faculty of Mechanics and  Mathematics, 
M.~V.~Lomonosov Moscow State University, Leninskie Gory 1, Moscow, 199991 Russia}

\author{Ol'ga Sipacheva}
\email[O. Sipacheva]{osipa@gmail.com}
\address{Department of General Topology and Geometry, Faculty of Mechanics and  Mathematics, 
M.~V.~Lomonosov Moscow State University, Leninskie Gory 1, Moscow, 199991 Russia}

\begin{abstract}
A topological space $X$ is \emph{$\mathbb R^{\omega_1}$-factorizable} if any continuous function $f\colon X\to 
\mathbb R^{\omega_1}$ factors through a continuous function from $X$ to a second-countable space. It is shown 
that a Tychonoff space $X$ is $\mathbb R^{\omega_1}$-factorizable if and only if $X\times D(\omega_1)$, where 
$D(\omega_1)$ is a discrete space of cardinality $\omega_1$, is $z$-embedded in the product $\beta X\times \beta 
D(\omega_1)$ of the Stone--\v Cech compactifications. It is also proved that $\mathbb 
R^{\omega_1}$-factorizability is hereditary and countably multiplicative, that any $\mathbb 
R^{\omega_1}$-factorizable space is hereditarily Lindel\"of and hereditarily separable, and that the existence of 
nonmetrizable $\mathbb R^{\omega_1}$-factorizable topological spaces and groups is independent of $\mathrm{ZFC}$: 
under $\mathrm{CH}$, all $\mathbb R^{\omega_1}$-factorizable spaces are second-countable, while under 
$\mathrm{MA}+\lnot\mathrm{CH}$, the countable Fr\'echet--Urysohn fan  is $\mathbb R^{\omega_1}$-factorizable. 
\end{abstract}



\maketitle

\section{Introduction}

A subset $S$ of a topological space $X$ is \emph{$z$-embedded} in $X$ if each zero set of $S$ is the trace on 
$S$ of a zero set of $X$. This notion has numerous applications in topology, analysis, and topological algebra 
(see, e.g., \cite{B-H, A-T}). In \cite{B-H} Blair and Hager studied the following property of the 
product $X\times Y$ of Tychonoff spaces: 
$$ 
\text{$X\times Y$ is $z$-embedded in $\beta X\times \beta Y$} \leqno(z) 
$$ 
(as usual, $\beta X$ and $\beta Y$ denote the Stone-\v Cech compactifications of $X$ and $Y$). They 
found conditions equivalent to (z) and considered applications of this property.

In \cite{R-S_2025} it was proved that a product $X\times Y$ has property (z) if and only if it is \emph{$\mathbb 
R$-factorizable}, that is, given any continuous function $f\colon X\times Y\to \mathbb R$, there exist separable 
metrizable spaces $X'$ and $Y'$ and continuous functions $g_X\colon X\to X'$, $g_Y\colon  Y\to Y'$, and $h\colon 
X'\times Y'\to \mathbb R$ such that $f= h\circ (g_X\times g_Y)$. The notion of $\mathbb R$-factorizability of 
products has applications in topological algebra, namely, in the theory of $\mathbb R$-factorizable groups 
\cite{R-S_2013, R_2024, R-S_2025} and $\mathbb R$-factorizable topological universal algebras 
\cite{R_2025}. 

It was shown in \cite{B-H} that if a product $X\times Y$ has property $(z)$ (or, equivalently, is 
$\mathbb R$-factorizable), then either $X$ or $Y$ is \emph{pseudo-$\aleph_1$-compact}, that is, contains no 
uncountable locally finite (or, equivalently, discrete) families of nonempty open sets. It is also known that if 
a product $X\times Y$ is Lindel\"of, then it is $\mathbb R$-factorizable \cite{R-S_2025}. The Sorgenfrey plane 
$S\times S$ does not have property (z) and hence is not $\mathbb R$-factorizable \cite{B-H}, but if $X$ 
is Lindel\"of and locally compact and $Y$ is pseudo-$\aleph_1$-compact, then $X\times Y$ is $\mathbb 
R$-factorizable \cite{B-H}. Moreover, a space $X$ is locally compact and second-countable if and only if 
$X\times Y$ is $z$-embedded in $\beta X\times \beta Y$  (or, equivalently, $\mathbb R$-factorizable) for any space $Y$~\cite{O}.

Given a cardinal $\kappa$, by $D(\kappa)$ we denote the discrete space of cardinality $\kappa$. In this paper 
we study spaces $X$ for which $X\times D(\omega_1)$ is $\mathbb R$-factorizable. They are interesting for several 
reasons, the main of which is that the $\mathbb R$-factorizability of $X\times D(\omega_1)$ is equivalent 
to the existence of a non-pseudo-$\aleph_1$-compact $Y$ for which $X\times Y$ is $\mathbb R$-factorizable 
\cite{R-S_2025}. It was also proved in \cite{R-S_2025} that $X\times D(\omega_1)$ is $\mathbb R$-factorizable if and only 
if, for any continuous function $f\colon  X\to \mathbb R^{\omega_1}$, there exists a separable metrizable space 
$Y$ and continuous functions $g\colon X\to Y$ and $h\colon Y\to \mathbb R^{\omega_1}$ such that $f=h\circ g$. For 
this reason, we refer to a space $X$ for which $X\times D(\omega_1)$ is $\mathbb R$-factorizable as an 
\emph{$\mathbb R^{\omega_1}$-factorizable} space. 

$\mathbb R^{\omega_1}$-Factorizable spaces have a number of interesting properties. For example, we prove in 
this paper that $\mathbb R^{\omega_1}$-factorizability is hereditary and countably multiplicative. In 
\cite{R-S_2025} we proved that any $\mathbb R^{\omega_1}$-factorizable space is hereditarily Lindel\"of and 
hereditarily separable and that any $\mathbb R^{\omega_1}$-factorizable space of weight at most $\omega_1$ is  
 second-countable. We also proved in \cite{R-S_2025} that $X\times D(2^\omega)$ is $\mathbb R$-factorizable 
if and only if $X$ is second-countable. Therefore, under $\mathrm{CH}$, all $\mathbb 
R^{\omega_1}$-factorizable spaces are metrizable (and hence second-countable). This suggests the question of 
whether the existence of nonmetrizable $\mathbb R^{\omega_1}$-factorizable spaces is consistent with 
$\mathrm{ZFC}$. In this paper we show that it is, which gives rise to the following problem. 

\begin{problemA}
\label{Problem_A}
What set-theoretic assumptions imply the existence of a nonmetrizable $\mathbb R$-factorizable space?
\end{problemA}

We give a partial solution of this problem. To this end, we use the following statement, which is implied by  
Propositions~\ref{Proposition_4.2} and~\ref{Proposition_4.5}, Corollary~\ref{Corollary_4.2}, and 
Theorem~\ref{Theorem_5.2} and allows us to seek 
nonmetrizable $\mathbb R$-factorizable spaces among spaces of simplest structure.

\begin{statementA}
\label{Statement_A}
The following assertions are equivalent:
\begin{enumerate}
\item[{\rm(1)}]
there exists a nonmetrizable $\mathbb R^{\omega_1}$-factorizable space;
\item[{\rm(2)}]
there exists a countable nonmetrizable $\mathbb R^{\omega_1}$-factorizable space with a single nonisolated point;
\item[{\rm(3)}]
there exists a countable nonmetrizable $\mathbb R^{\omega_1}$-factorizable space of weight $\omega_2$ with a 
single nonisolated point; 
\item[{\rm(4)}] 
there exists a nonmetrizable $\mathbb R^{\omega_1}$-factorizable 
topological group. 
\end{enumerate} 
\end{statementA}

If $X$ is an $\mathbb R^{\omega_1}$-factorizable space and $x\in X$, then the neighborhood filter $\mathscr U_x$ 
of $x$ (and the filter $\mathring {\mathscr U}_x$ of punctured neighborhoods of $x$ if $x$ 
is nonisolated) has the following property:
given any $U_\alpha\in \mathscr U$, $\alpha<\omega_1$, there exist $V_n\in \mathscr U$, $n<\omega$, 
such that every $U_\alpha$ contains $V_{n_\alpha}$ for some $n_\alpha<\omega$ 
(see Proposition~\ref{Proposition_4.1}).
We refer to such filters as \emph{$A$-filters}.

Given a free filter $\mathscr F$ on a set $X$, by $X_{\mathscr F}$ we denote the space $X\cup\{\infty_X\}$, where 
$\infty_X\notin X$, in which all points of $X$ are isolated and $\mathscr F$ is the filter of punctured 
neighborhoods of~$\infty_X$. Any space with a single nonisolated point is homeomorphic to a space of the form 
$X_{\mathscr F}$. Clearly, $X_{\mathscr F}$ is metrizable if and only if the filter $\mathscr F$ has a countable 
base. It turns out that, for countable $X$, $X_{\mathscr F}$ is $\mathbb R^{\omega_1}$-factorizable if and only 
if $\mathscr F$ is an $A$-filter. Our results imply the following statement (see Proposition~\ref{Proposition_4.1} 
and Corollary~\ref{Corollary_4.2}).

\begin{statementA}
\label{Statement_B}
Let $\mathscr F$ be a free filter on $\omega$. The space $\omega_{\mathscr F}$ is 
$\mathbb R^{\omega_1}$-factorizable if and only if $\mathscr F$ is an $A$-filter.
The space $\omega_{\mathscr F}$ is nonmetrizable if and only if the filter $\mathscr F$ has no 
countable base. 
\end{statementA}

It follows from Statements~\ref{Statement_A} and~\ref{Statement_B} that Problem~\ref{Problem_A} 
is equivalent to the following one.

\begin{problemA}
What set-theoretic assumptions imply the existence of $A$-filters on $\omega$ without a countable base? 
\end{problemA}

We prove that one of such assumptions is $\mathfrak b>\omega_1$ (see Proposition~\ref{Proposition_4.4}), although the 
existence of $A$-filters on $\omega$ without a countable base is consistent with $\mathfrak b=\omega_1$ (see 
Corollary~\ref{Corollary_4.3}). The key assumption is the existence of a strictly $\subset^*$-descending 
$\omega_2$-sequence of subsets of $\omega$. More precisely, we prove the following 
statement  (see Theorem~\ref{Theorem_4.1}, Corollaries~\ref{Corollary_4.2} and~\ref{Corollary_4.3}, 
and Proposition~\ref{Proposition_4.4}). 

\begin{statementA}
\begin{enumerate}
\item[{\rm(1)}]
The existence of a strictly $\subset^*$-descending 
$\omega_2$-sequence of subsets of $\omega$ implies that of a nonmetrizable $\mathbb R^{\omega_1}$-factorizable 
space. 
\item[{\rm(2)}] 
The assumption $\mathfrak b>\omega_1$ implies the existence of a strictly $\subset^*$-descending 
$\omega_2$-sequence of subsets of $\omega$ and hence the existence of a nonmetrizable $\mathbb 
R^{\omega_1}$-factorizable space. Moreover, $\mathfrak b>\omega_1$ if and only if the Fr\'echet--Urysohn 
fan $V(\omega)$ is $\mathbb R^{\omega_1}$-factorizable. 
\item[{\rm(3)}] 
The existence of a strictly $\subset^*$-descending 
$\omega_2$-sequence of subsets of $\omega$ is consistent with $\omega_1=\mathfrak b<2^\omega=\omega_2$. 
\end{enumerate}
\end{statementA}

We also consider $\mathbb R^{\omega_1}$-factorizable topological groups. We prove that a topological group $G$ is 
$\mathbb R^{\omega_1}$-factorizable if and only if any continuous function $f\colon G\to \mathbb R^{\omega_1}$ 
factors through a homomorphism to a second-countable group, i.e., there exists a second-countable 
topological group $H$, a continuous homomorphism $h\colon G \to H$, and a continuous function $g\colon H\to 
\mathbb R^{\omega_1}$ for which $f = g \circ h$. In particular, any $\mathbb R^{\omega_1}$-factorizable group is 
$\mathbb R$-factorizable. We prove that the existence of a nonmetrizable $\mathbb R^{\omega_1}$-factorizable 
topological group is equivalent to that of an $A$-filter without a countable base on~$\omega$. 

The referee kindly informed us of the joint  paper \cite{HPTZ} by Wei He, Dekui Peng, Mikhail Tkachenko, and 
Heng Zhang, in which $\tau$-fine topological groups were introduced and studied. An $\omega_1$-fine group is 
precisely a group in which the filter of neighborhoods of the identity element is an $A$-filter. Our results and 
those of \cite{HPTZ} complement each other very well (see, e.g., our Theorem~\ref{Theorem_5.3}). 

The concluding section of the paper is devoted to open problems.

\section{Preliminaries} 

For simplicity, we assume all topological spaces under consideration to be Tychonoff and often refer to 
them simply as spaces.

\begin{definition}[\cite{R-S_2013, R_2024}]
\label{Definition_2.1}
Given topological spaces $X$, $Y$, and $Z$, we say that a function $f\colon X\times Y\to Z$ 
is \emph{$\mathbb R$-factorizable} if it \emph{factors through a product} of continuous functions to 
second-countable spaces, that is, if there exist second-countable spaces $X'$ and $Y'$ and continuous 
functions $g\colon X\to X'$, $h\colon Y\to Y'$, and $\varphi\colon X'\times Y'\to Z$ such that 
$f=\varphi\circ(g\times h)$, i.e., the following diagram 
is commutative: 
$$
\begin{tikzcd} 
&X\times Y\arrow[r, "f"] \arrow[d,shift left=12pt]\arrow[d,shift right=12pt,"\rlap{$g\,\times\, h$}"] & Z\\
&X'\times Y'\arrow[ru, "\varphi"']& 
\end{tikzcd} 
$$
We say that a product $X\times Y$ is \emph{$\mathbb R$-factorizable} 
(or \emph{topologically $\mathbb R$-fac\-tor\-iz\-able}, when there is a danger of confusion) 
if any continuous function $f\colon X\times Y \to \mathbb R$ is $\mathbb R$-factorizable. 
\end{definition} 

\begin{definition}
Given a cardinal $\kappa$, a topological space $X$ is said to be \emph{pseudo-$\kappa$-compact} if the 
cardinality of any locally finite family of open sets in $X$ is less than~$\kappa$. 
\end{definition}

\begin{remark}[{see, e.g., \cite{R-S_2025}}]
It is easy to see that a Tychonoff space $X$ is pseudo-$\kappa$-compact if and only if the cardinality of 
any discrete family of cozero sets in $X$ is less than~$\kappa$. 
\end{remark}

Given a cardinal $\kappa$, we use $D(\kappa)$ to denote a discrete space of cardinality~$\kappa$. 
We denote the 
weight and character of a space $X$ by~$w(X)$ and $\chi(X)$, respectively. 
If $\{X_\iota: \iota\in I\}$ is a family of sets or spaces,
$J\subset S\subset I$ are arbitrary subsets of the index set, and $\nu\in I$, then by 
$\pi_J^S$, $\pi_J$, and $\pi_\nu$ 
we denote the natural projection maps 
$\prod_{\iota\in S}X_\alpha \to \prod_{\iota\in J}X_\iota$,
$\prod_{\iota\in I}X_\alpha \to \prod_{\iota\in J}X_\iota$, and $\prod_{\iota\in I}X_\iota \to X_\nu$, 
respectively. 
For the cardinality of a set $X$ we use the notation $|X|$. If 
$Y$ is a subset of a topological space 
$X$, then by $\overline Y$ we denote the closure of $Y$ in $X$. We use the symbol $\oplus$ for topological sum. 
By $\omega^\omega$ We denote the set of all functions $\omega\to \omega$ and by $[\omega]^\omega$, the set of all 
infinite subsets of~$\omega$. By a \emph{cosmic} space we mean a space with a countable network.

The following theorem collects some of the results obtained in~\cite{R-S_2025} (see 
Propositions~2.2--2.4 of~\cite{R-S_2025}). 

\begin{theoremA}[\cite{R-S_2025}]
\label{Theorem_A}
For a Tychonoff space $X$, the following conditions are equivalent:
\begin{enumerate}
\item[\rm(1)]
the product $X\times D(\kappa)$ is $\mathbb R$-factorizable; 
\item[\rm(2)]
there exists a non-pseudo-$\kappa$-compact space $Y$ for which $X\times Y$ is $\mathbb R$-factorizable; 
\item[\rm(3)] 
every continuous function $f\colon X\to Y$ to a Tychonoff space $Y$ with $w(Y)\leq\kappa$ factors 
through a continuous function to a second-countable space.
\end{enumerate} 
\end{theoremA}

\section{$\mathbb R^{\omega_1}$-Factorizable Spaces}

\begin{definition}
We say that a topological space $X$ is \emph{$\mathbb R^{\omega_1}$-factorizable} if any continuous function  
$f\colon X\to \mathbb R^{\omega_1}$ factors through a continuous function to a second-countable space, 
that is, given any continuous function $f\colon X\to \mathbb R^{\omega_1}$, there exists 
a second-countable space $Y$, a continuous function $h\colon X \to Y$, and a continuous function 
$g\colon Y\to \mathbb R^{\omega_1}$ for which $f = g \circ h$, i.e., the following diagram is commutative:
$$
\begin{tikzcd} 
&X \arrow[r, "f"]  \arrow[d, "h"']   & \mathbb R^{\omega_1} \\ 
&Y \arrow[ru, "g"'] &
\end{tikzcd} 
$$
\end{definition}

\begin{proposition}
A space $X$ is $\mathbb R^{\omega_1}$-factorizable if and only if 
the product $X\times D(\omega_1)$ is $\mathbb R$-factorizable. 
\end{proposition}

\begin{proof}
Since any Tychonoff space $Y$ with $w(Y)\le \omega_1$ embeds in $\mathbb R^{\omega_1}$, it follows that a 
space $X$ is $\mathbb R^{\omega_1}$-factorizable if and only if any continuous function from $X$ to a 
Tychonoff space of weight at most $\omega_1$ factors through a continuous function to a second-countable space. 
It remains to apply Theorem~\ref{Theorem_A}\,(3).
\end{proof}

The following theorem is an immediate corollary of Theorem~\ref{Theorem_A}. 

\begin{theorem}
\label{Theorem_3.1}
A space $X$ is $\mathbb R^{\omega_1}$-factorizable if and only if 
there exists a non-pseudo-$\aleph_1$-compact space $Y$ for which $X\times Y$ is $\mathbb R$-factorizable.
\end{theorem}

The following criterion plays a key role in the study of $\mathbb R^{\omega_1}$-factorizable spaces.

\begin{proposition}
\label{Proposition_3.2}
Let $\{M_\iota:\iota\in I\}$ a family of second-countable spaces and $X\subset \prod_{\iota\in I}M_\iota$.
Then following conditions are equivalent:
\begin{enumerate}
\item[\rm (1)]
$X$ is $\mathbb R^{\omega_1}$-factorizable;
\item[\rm (2)]
for any $J\subset I$ with $|J|\le\omega_1$, there exists a second-countable space 
$M$ and continuous maps $g\colon X\to M$ and $h\colon M\to \prod_{\iota\in J}M_\iota$ 
for which $\pi_J|_{X} = h\circ g$; 
\item[\rm (3)]
for any $J\subset I$ with $|J|\le\omega_1$, there exists  $S\subset I$ for which $J\subset S$ and $\pi_S(X)$ is second-countable;
\item[\rm (4)]
for any $J\subset I$ with $|J|\le\omega_1$, there exists at most countable $Q\subset I$ 
and a continuous map $\rho\colon \pi_Q(X)\to \pi_J(X)$
for which 
$\pi_J|_X = \rho\circ \pi_Q|_X$.
\end{enumerate}
\end{proposition}
\begin{proof}
First, we prove two claims.

\begin{claim-1} 
If condition $(2)$ holds, then $X$ is Lindel\"of.
\end{claim-1}

\begin{proof}
Suppose that there 
exists an open cover $\mathscr U$ of $X$ by open subsets of $\prod_{\iota\in I}M_\iota$ which has no countable 
subcover. We can assume without loss of generality 
that $\mathscr U=\{U(\alpha): \alpha\in \kappa\}$, where $\kappa$ is an uncountable cardinal and 
each $U(\alpha)$ is a canonical open set in 
$\prod_{\iota\in I}M_\iota$, that is, $U(\alpha)$ is the product $\prod_{\iota\in I}V_\iota(\alpha)$ of 
open subsets $V_\iota(\alpha)$ of $M_\iota$ and $V_\iota(\alpha)=X_\iota$ for all $\iota\in I\setminus 
F_\iota$, where $F_\iota$ is a finite subset of $I$. We can also assume that, for each $\alpha<\omega_1$, 
$X\cap U(\alpha)\not\subset X\cap \bigcup_{\beta<\alpha} U(\beta)$. Let $J=\bigcup_{\alpha<\omega_1}F_\alpha$.
Then $|J|\le \omega_1$ and by assumption there exists a second-countable space 
$M$ and continuous maps $g\colon X\to M$ and $h\colon M\to \prod_{\iota\in J}M_\iota$ 
for which $\pi_J|_{X} = h\circ g$. 
Therefore, $\pi_J(X)$ is cosmic, being a continuous image of a second-countable space $M$. 
For $\gamma<\omega_1$, we set $V(\gamma)=\pi_J(U(\gamma))=\prod_{\iota\in J}V_\iota(\gamma)$.
Since $F_\gamma\subset J$, it follows that $U(\gamma)=\pi_J^{-1}(V(\gamma))$.
Hence $\pi_J(X)\cap V(\alpha)\not\subset \pi_J(X)\cap  \bigcup_{\beta<\alpha}  V(\beta)$. 
Thus, $\{\pi_J(X)\cap V(\alpha): \alpha<\omega_1\}$ is an uncountable open cover with 
no countable subcover of the cosmic space 
$\pi_J(X)\cap \bigcup\{V(\alpha): \alpha<\omega_1\}$. 
This contradiction proves that $X$ is Lindel\"of. 
\end{proof}

\begin{claim-2}
If $X$ is Lindel\"of, $Y$ is any infinite space, and $f\colon X\to Y$ is a continuous map, 
then there exists a $J\subset I$ with $|J|\leq w(Y)$ and a continuous map $h\colon \pi_J(X)\to Y$ 
for which $f=h\circ \pi_J|_X$.
\end{claim-2}

\begin{proof}
In view of Tychonoff's embedding theorem we can assume that $Y\subset {\mathbb R}^\kappa$ for $\kappa=w(Y)$.
According to Lemma~8.1.4 of \cite{A-T}, for
every function $\pi_\alpha\circ f\colon X\to \mathbb R$, $\alpha<\kappa$, there exists a countable set 
$J_\alpha\subset I$ and a continuous function $\varphi_\alpha\colon \pi_{J_\alpha}(X)
\to \mathbb R$ such that $\pi_\alpha\circ f = \varphi_\alpha\circ \pi_{J_\alpha}|_{X}$. 
Let $J=\bigcup_{\alpha< \kappa} J_\alpha$. Then $|J|\le \kappa = w(Y)$. 
For  each $\alpha<\kappa$, let $h_\alpha\colon \pi_J(X)\to \mathbb R$ 
defined by $h_\alpha(x)=\varphi_\alpha(\pi^J_{J_\alpha}(x))$ for $x\in \pi_J(X)$.
Then $h_\alpha=\pi_\alpha\circ f\circ \pi_J|_X$ 
for $\alpha<\kappa$, the map $h=\bigdiag_{\alpha<\kappa}h\colon \pi_J(X)\to {\mathbb R}^\kappa$ 
is continuous, and $f=h\circ \pi_J|_X$.
\end{proof}

We proceed to the proof of the proposition. 	       
The implication (1)\,$\Rightarrow$\,(2) holds by the definition of $\mathbb R^{\omega_1}$-factorizable spaces. 

Let us show that (2)\,$\Rightarrow$\,(4). Let $J$ be a subset of $I$ with $|J|\le \omega_1$.
By (2) there exists a second-countable space $M$ and continuous maps $g\colon X\to M$ and 
$h\colon M\to \pi_J(X)\subset \prod_{\iota\in J}M_\iota$ for which $\pi_J|_{X} = h\circ g$.
According to Claim~1, $X$ is Lindel\"of, and Claim~2 implies the existence of an at most countable set 
$Q\subset I$ and a continuous map $\psi\colon \pi_Q(X)\to M$ such that $g=\psi\circ \pi_Q|_{X}$.
We have $\pi_J|_{X} = h\circ \psi\circ \pi_Q|_{X}$. It remains to set $\rho=h\circ \psi$.

Let us prove the implication (4)\,$\Rightarrow\,$(3). Let $J$ be a subset of $I$ with $|J|\le \omega_1$.
By (4) there exists an at most countable set $Q\subset I$ and a continuous map $\rho\colon \pi_Q(X)\to \pi_J(X)$
for which $\pi_J|_X = \rho\circ \pi_Q|_X$.
Let $S=J\cup Q$. Consider the map $h\colon \pi_Q(X)\to \pi_S(X)$ defined by
\[
h(x)(\iota) = \begin{cases}
x(\iota) &\text{if }\iota\in Q\setminus J,
\\
\rho(x)(\iota) &\text{if }\iota\in J
\end{cases}
\]
for $x\in \pi_Q(X)$.
This map is continuous, being the diagonal of the identity map $\pi_{Q\setminus J}\to \pi_{Q\setminus J}$ 
and $\rho$, and it is inverse to $\pi^S_Q|_{\pi_S(X)}\colon \pi_S(X)\to \pi_Q(X)$. Therefore, 
$\pi^S_Q|_{\pi_S(X)}$ is a homeomorphism, so that the spaces $\pi_S(X)$ and $\pi_Q(X)$ are homeomorphic. 
It follows that $\pi_S(X)$ is second-countable.

Let us prove (3)\,$\Rightarrow$\,(2). Given $J\subset I$ with $|J|\le\omega_1$, 
take  $S\subset I$ for which $J\subset S$ and $M=\pi_S(X)$ is second-countable and 
let $g=\pi_S|_X$ and $h=\pi^S_J|_M$. Then $\pi_J|_{X} = h\circ g$.

Finally, we prove that (2)\,$\Rightarrow$\,(1). 
According to Claim 1, $X$ is Lindel\"of.
Let $f\colon X\to \mathbb R^{\omega_1}$ be a continuous function.
Claim~2 implies the existence of a set $J\subset I$, $|J|\leq \omega_1$, and a continuous map 
$\psi\colon \pi_J(X)\to \mathbb R^{\omega_1}$ for which $f=\psi\circ \pi_J|_X$.
By condition (2) 
there exists a second-countable space 
$M$ and continuous maps $h\colon X\to M$ and $\varphi\colon M\to \pi_J(X)$ 
for which $\pi_J|_{X} = \varphi\circ h$. 
Let $g=\psi\circ \varphi \colon M\to \mathbb R^{\omega_1}$.
Then $f = g\circ h$.
\end{proof}

The following statement is essentially Proposition~2.5 of \cite{R-S_2025}, but its proof given in 
\cite{R-S_2025} contains an easy-to-fill gap\footnote{In the proof of Proposition~2.5, it is claimed that 
$\bigdiag_{\alpha<\omega_1}f_\alpha^{\times\omega}(F_\beta)$ is closed in 
$\prod_{\alpha<\omega_1}M_\alpha^\omega$. In fact, 
it is closed in $\bigdiag_{\alpha<\omega_1}f_\alpha^{\times\omega}(X^\omega)$. 
This does not affect further argument.}. Below we 
suggest a simpler alternative proof. 

\begin{corollary}
\label{Corollary_3.1}
If a space $X$ is $\mathbb R^{\omega_1}$-factorizable, then so is $X^\omega$.
\end{corollary}

\begin{proof}
Let $X$ be an $\mathbb R^{\omega_1}$-factorizable space; we assume it to be infinite. 
By Tychonoff's embedding theorem, 
there is an embedding $i\colon X\to \mathbb R^\kappa$, where $\kappa=w(X)$. Consider the product map 
$i^{\times \omega}\colon X^\omega\to \mathbb R^{\kappa\times\omega}$; this is an embedding. 

According to Proposition~\ref{Proposition_3.2}\,(3), to prove that $X^\omega$ is 
$\mathbb R^{\omega_1}$-factorizable, it suffices to show that, for any $J\subset \kappa\times \omega$ with 
$|J|\le \omega_1$, there exists an $S\subset \kappa\times \omega$ for which $J\subset S$ and 
$\pi_S^{\kappa\times \omega}(X^\omega)$ is second-countable.

There exists an $R\subset \kappa$ with $|R|\leq\omega_1$ such that $J\subset R\times \omega$.
Proposition~\ref{Proposition_3.2}(3) implies the existence of a $T\subset \kappa$ with $|T|\leq\omega_1$ 
for which $R\subset T$ and $\pi_T^\kappa(X)$ is second-countable.
Let $S=T\times \omega$. Then $\pi_S^{\kappa\times \omega}(X^\omega)=\pi_T^\kappa(X)^\omega$ is second-countable.
\end{proof}

The following theorem is an immediate consequence of Proposition~2.6 of~\cite{R-S_2025}.

\begin{theorem}
\label{Theorem_3.2}
Let $X$ be an $\mathbb R^{\omega_1}$-factorizable space. Then
\begin{enumerate}
\item[\rm(1)]
$X^\omega$ is hereditarily Lindel\"of and hereditarily separable; 
\item[\rm(2)]
if $w(X)\leq\omega_1$, then $X$ is second-countable.
\end{enumerate}
\end{theorem}

\begin{theorem}
\label{Theorem_3.3}
The following assertions hold.
\begin{enumerate}
\item[\rm(1)]
Any subspace of an $\mathbb R^{\omega_1}$-factorizable space is $\mathbb R^{\omega_1}$-factorizable. 
\item[\rm(2)]
If $X_n$, $n\in \omega$, are $\mathbb R^{\omega_1}$-factorizable spaces, then 
$\bigoplus_{n\in \omega} X_n$ is $\mathbb R^{\omega_1}$-factorizable.
\item[\rm(3)]
If $X_n$, $n\in \omega$, are $\mathbb R^{\omega_1}$-factorizable spaces, then 
$\prod_{n\in \omega} X_n$ is $\mathbb R^{\omega_1}$-fac\-tor\-iz\-able.
\end{enumerate}
\end{theorem}

\begin{proof}
Let us prove~(1). 
Let $X$ be an $\mathbb R^{\omega_1}$-factorizable space, and let $Y \subset X$. 

We can assume that $X$ is a subspace of $\mathbb R^{I}$ for some set $I$.
According to Proposition~\ref{Proposition_3.2}\,(3), to prove that $Y$ is $\mathbb R^{\omega_1}$-factorizable, 
it suffices to show that, for any $J\subset I$ with $|J|\le \omega_1$, there exists an $S\subset I$ 
for which $J\subset S$ and $\pi_S(Y)$ is second-countable.
By the same Proposition~\ref{Proposition_3.2}\,(3) there exists an $S\subset I$ for which $J\subset S$ 
and $\pi_S(X)$ is second-countable. Clearly, $\pi_S(Y)$ is second-countable.

Assertion~(2) is almost obvious. Indeed, given a continuous function $f\colon \bigoplus_{n\in \omega}X_n\to 
\mathbb R^{\omega_1}$, consider the restrictions $f_n=f|_{X_n}\colon X_n\to \mathbb R^{\omega_1}$. By assumption 
for each $n\in \omega$ there exists a second-countable space $M_n$ and continuous functions $h_n\colon X_n\to 
M_n$ and $g_n\colon M_n \to \mathbb R^{\omega_1}$ such that $f_n=g_n\circ h_n$. Let $h\colon \bigoplus_{n\in 
\omega}X_n\to \bigoplus_{n\in \omega} M_n$ and $g\colon \bigoplus_{n\in \omega} M_n\to \mathbb R^{\omega_1}$ be 
functions uniquely determined by the conditions $h|_{X_n}=h_n$ and $g|_{M_n}=g_n$, respectively. Clearly, $h$ and 
$g$ are continuous, $\bigoplus_{n\in \omega} M_n$ is second-countable, and $f=g\circ h$. 

Assertion~(3) follows from (1), (2), Corollary~\ref{Corollary_3.1}, and the 
inclusion $\prod_{n\in \omega}X_n\subset 
\bigl(\bigoplus_{n\in \omega}X_n\bigr)^\omega$.
\end{proof}

\begin{remark}
Let $X$ be a countable space of weight $\omega_1$. Then $X$ is not $\mathbb 
R^{\omega_1}$-factorizable by Theorem \ref{Theorem_3.2}\,(2). The space $X$ is a continuous image of a countable discrete space, which is  
$\mathbb R^{\omega_1}$-factorizable. Therefore, $\mathbb R^{\omega_1}$-factorizability is not preserved by 
continuous maps. 
\end{remark}

\section{$A$-Filters}

\begin{definition}
We say that a filter $\mathscr F$ on an arbitrary set is an \emph{$A$-filter} if, given any 
$F_\alpha\in \mathscr F$, $\alpha<\omega_1$, there exist $G_n\in \mathscr F$, $n<\omega$, such that 
each $F_\alpha$ contains $G_{n_\alpha}$ for some $n_\alpha<\omega$. 
\end{definition}

\begin{proposition} 
\label{Proposition_4.1}
The neighborhood filter $\mathscr U_x$ of any point $x$ in an $\mathbb R^{\omega_1}$-factorizable space~$X$ 
is an $A$-filter.
\end{proposition}

\begin{proof}
Let $U_\alpha\in \mathscr U_x$, $\alpha<\omega_1$.
For each $\alpha<\omega_1$, there exists a  continuous function $f_\alpha$ on $X$ such that 
$f_\alpha(x)=0$ and $f_\alpha(X\setminus U_\alpha)\subset\{1\}$.
Consider the function $f=\diag_{\alpha<\omega_1}\colon X\to R^{\omega_1}$. The 
$\mathbb R^{\omega_1}$-factorizability of $X$ implies the existence of a second-countable space $Y$ and 
continuous functions $h\colon X \to Y$ and $g\colon Y\to \mathbb R^{\omega_1}$ 
for which $f = g \circ h$.
Let $\{W_n:n\in\omega\}$ be a local base at $h(x)\in Y$, and let $V_n=h^{-1}(W_n)$ for $n\in\omega$.
Then every $U_\alpha$ contains some~$V_n$.
\end{proof}

\begin{proposition}
\label{Proposition_4.2}
A countable space $X$ 
is $\mathbb R^{\omega_1}$-factorizable if and only if 
the neighborhood filter of every point $x\in X$ is an~$A$-filter. 
\end{proposition}

\begin{proof}
Let us prove the `if' part. Suppose that the neighborhood filter $\mathscr U_x$ of every
point $x\in X$ is an $A$-filter. 
Let $Z$ be a space with a base $\mathscr{B}$, $|\mathscr{B}|\leq\omega_1$, and let $f\colon X\to Z$ be a
continuous function. We can assume without loss of generality that $f$ is surjective. 
We set $\mathscr{W}=\{ f^{-1}(U): U\in \mathscr{B}\}$ and 
$\mathscr{W}_x=\{W\in \mathscr{W}: x\in W\}$ for $x\in X$.
The space $X$ is zero-dimensional and has countable pseudocharacter, being countable.
For each $x\in X$, there exists a family $\{V_{x,n}:n\in\omega\}\subset \mathscr U_x$ of clopen subsets of $X$ 
such that $\bigcap_{n\in \omega}V_{x,n}=\{x\}$ and every $W\in \mathscr{W}_x$ contains $V_{x,n}\subset W$ 
for some $n\in\omega$, because $\mathscr U_x$ is an $A$-filter. 
We set $\mathscr{V}=\{V_{x,n}: x\in X, n\in\omega\}$ and 
$h=\bigdiag_{V\in \mathscr{V}} \chi_V\colon X\to \{0,1\}^{\mathscr{V}}$, where $\chi_V$ denotes 
the characteristic function 
of $V\in \mathscr{V}$. 
Then $h$ is injective, $Y=h(X)$ is second-countable, and the map  
$g\colon Y\to Z$ defined by $g(y)=f(h^{-1}(y))$ is continuous. Indeed, for any $y\in Y$ and $U\ni g(y)$, 
$U\in \mathscr B$, there exists an $n\in \omega$ such that $V_{h^{-1}(y),n}\subset W=f^{-1}(U)$. 
We have $g(z)\in U$ for all $z\in Y$ such that $z(V_{h^{-1}(y),n})=1$. Therefore, $g$ is continuous at 
every $y\in Y$. Clearly, $f=g\circ h$.

The reverse implication follows from Proposition~\ref{Proposition_4.1}.
\end{proof}

Given a free filter $\mathscr F$ on a set $X$, by $X_{\mathscr F}$ we denote the space $X\cup\{\infty_X\}$, where 
$\infty_X\notin X$, in which all points of $X$ are isolated and $\mathscr F$ is the filter of punctured 
neighborhoods of~$\infty_X$. Proposition~\ref{Proposition_4.2} has the following immediate corollary. 

\begin{corollary}
\label{Corollary_4.1}
 A free filter $\mathscr F$  on a countable set $X$ is an $A$-filter if and only if 
the space $X_{\mathscr F}$ is $\mathbb 
R^{\omega_1}$-factorizable.
\end{corollary}

Thus, to construct an example of a nonmetrizable $\mathbb R^{\omega_1}$-factorizable space, it suffices to 
construct a free $A$-filter on a countable set which has no countable base. To this end, we use the preorder 
$\subset^*$ on $[\omega]^\omega$ defined by setting $A\subset^* B$ if $A\setminus B$ is finite. Given a cardinal 
$\kappa$, a family $\{B_\alpha: \alpha<\kappa\} \subset [\omega]^\omega$ is called a \emph{$\subset^*$-descending 
$\kappa$-sequence} if $B_\beta\subset^* B_\alpha$ whenever $\alpha<\beta<\kappa$ and a \emph{strictly 
$\subset^*$-descending $\kappa$-sequence} if $\mathscr B$ is $\subset^*$-descending and $B_\alpha\setminus 
B_\beta$ is infinite whenever $\alpha<\beta<\kappa$. 

\begin{lemma}
\label{Lemma_4.1}
Let $\kappa\ge \omega_2$ be a regular cardinal, and let $\mathscr B=\{B_\alpha: \alpha<\kappa\}$ 
be a $\subset^*$-descending $\kappa$-sequence in $[\omega]^\omega$. Then the filter $\mathscr F$ on 
$\omega$ generated by the subbase 
$$ 
\mathscr B' = \{\omega \setminus \{n\}: n\in\omega\}\cup \mathscr B
$$
is an $A$-filter. If, in addition, $\mathscr B$ is strictly $\subset^*$-descending, then $\mathscr F$ has no 
countable base. 
\end{lemma}

\begin{proof}
Consider any family $\{F_\alpha: 
\alpha<\omega_1\}\subset \mathscr F$. Note that, for every $\alpha<\omega_1$, there exists a $\beta<\kappa$ such 
that $B_\beta\subset^* F_\alpha$. Since $\kappa$ is regular and $\kappa>\omega_1$, it follows that there 
exists an $\alpha_0<\kappa$ such that $B_{\alpha_0}\subset^* F_\alpha$ for all $\alpha<\kappa$. Setting 
$B'_n=B_{\alpha_0}\setminus \{0,\dots, n\}$ for $n<\omega$, we obtain the countable set $\{B'_n: n<\omega\}$ of 
elements of $\mathscr F$ satisfying the condition in the definition of an $A$-filter. Thus, $\mathscr F$ is an 
$A$-filter. The second assertion is obvious. 
\end{proof}

Note that under $\mathrm{CH}$ any filter on $\omega$ contains at most $\omega_1$ elements. Therefore, for any 
$A$-filter $\mathscr F$ on a countable set $X$, we have $w(X_{\mathscr F})\le\omega_1$; it follows from 
Theorem~\ref{Theorem_3.2}\,(2) that any such filter has a countable base. This observation, together with Lemma~\ref{Lemma_4.1},  
implies the following theorem.

\begin{theorem}
\label{Theorem_4.1}
\begin{enumerate}
\item[\rm(1)]
If $[\omega]^\omega$ contains a strictly $\subset^*$-de\-scend\-ing $\omega_2$-sequence, then there exists an 
$A$-filter on $\omega$ which has no countable base. 
\item[\rm(2)] 
If $2^\omega=\omega_1$, then any $A$-filter on $\omega$ has a countable base. 
\end{enumerate} 
\end{theorem}

Recall that $\le^*$ denotes the preorder on $\omega^\omega$ defined by setting $f\le^* g$ if $f(n)\le g(n)$ for 
all but finitely many $n \in \omega$. A subset of $\omega^\omega$ is said to be \emph{unbounded} if it is 
unbounded with respect to this preorder. The cardinal $\mathfrak b$ is defined by
$$
\mathfrak b = \min\{|B|: \text{$B$ is an unbounded subset of $\omega^\omega$}\}.
$$

\begin{proposition}
There exists a strictly $\subset^*$-de\-scend\-ing $\mathfrak b$-sequence in~$[\omega]^\omega$.
\end{proposition}

\begin{proof}
Using the fact that any subset of $\omega^\omega$ of cardinality less than $\mathfrak b$ is bounded, we can 
easily construct by recursion a set of functions $f_\alpha\in \omega^\omega$, $\alpha<\mathfrak b$, such that 
$f_\alpha(n) <^* f_\beta(n)$ for all but finitely many $n\in \omega$ whenever $\alpha<\beta<\mathfrak b$. Setting $B_\alpha=\{(m,n): m\in \omega,\ n\ge 
f_\alpha(m)\}$ and identifying $\omega\times \omega$ with $\omega$, we obtain the desired strictly 
$\subset^*$-de\-scend\-ing $\mathfrak b$-sequence. 
\end{proof}

\begin{corollary}
\label{Corollary_4.2}
If $\mathfrak b>\omega_1$, then there exists a countable nonmetrizable $\mathbb R^{\omega_1}$-factorizable space.
\end{corollary}

In \cite{B-M} Br\"auninger and Mildenberger constructed a ZFC model in which $\omega_1=\mathfrak 
b<2^\omega=\omega_2$ and there exists a simple $P_{\aleph_2}$-point in $\omega^*$, that is, an ultrafilter 
$\mathscr U$ on $\omega$ which has a base being a strictly $\subset^*$-descending $\omega_2$-sequence. 
This implies the following statement in view of Theorem~\ref{Theorem_4.1}\,(1) and Corollary~\ref{Corollary_4.2}.

\begin{corollary}
\label{Corollary_4.3}
Both assumptions $\mathfrak b=\omega_1$ and $\mathfrak b>\omega_1$ are consistent with the existence of a 
countable nonmetrizable $\mathbb R^{\omega_1}$-factorizable space.
\end{corollary}

Let $V(\omega)$ be the \emph{Fr\'echet--Urysohn fan}, that is, the image of the sum of countably many copies of 
the convergent sequence $\{\frac 1n: n\in \mathbb N\}\cup \{0\}$ under the quotient map contracting the set of 
limit points of all these sequences to a single point~$O$. 

\begin{proposition}
\label{Proposition_4.4}
The following conditions are equivalent:
\begin{enumerate}
\item[\rm(1)]
$\mathfrak b>\omega_1$\textup;
\item[\rm(2)]
the filter of neighborhoods of $O$ in $V(\omega)$ is an $A$-filter and has no countable base;
\item[\rm(3)]
$V(\omega)$ is a nonmetrizable $\mathbb R^{\omega_1}$-factorizable space.
\end{enumerate}
\end{proposition}

\begin{proof}
The fan $V(\omega)$ can be represented as 
$$
V(\omega)=\{O\}\cup \{(n, 1/k): n\in \omega,\ k\in \mathbb N\}.
$$
The filter $\mathscr F$ of neighborhoods of $O$ has the base formed by the sets 
$$
F_f=\{(n, 1/k)\cup\{O\}: n\in \omega,\ k\ge f(n)\},\qquad f\in \omega^\omega. 
$$

Let us prove the implication (1)\,$\Rightarrow$\,(2).
Consider a family $\{F_{f_\alpha}: \alpha <\omega_1\}$ of elements of $\mathscr F$. 
Since $\mathfrak b>\omega_1$, it follows that there exists a function $g\in 
\omega^\omega$ such that $f_\alpha \le^* g$ for each $\alpha<\omega_1$. For every pair $(l,m)\in \omega\times 
\omega$, we set 
$$
g_{(l,m)}(n)=
\begin{cases}
l&\text{if $n\le m$},\\
g(n)&\text{otherwise}.
\end{cases}
$$  
For each $\alpha<\omega_1$, there exists an $(l,m)\in \omega$ such that $f_\alpha(n)\le g_{(l,m)}(n)$ for 
all $n\in \omega$ (because $f_\alpha\le^* g$) and hence $F_{g_{(l,m)}}\subset F_{f_\alpha}$. 

The filter $\mathscr F$ has no countable base, because the fan $V(\omega)$ is not first-countable.

The equivalence (2)\,$\Leftrightarrow$\,(3) follows from Proposition~\ref{Proposition_4.2}. 

Finally, we prove the implication (2)\,$\Rightarrow$\,(1) by contradiction. Suppose that (2) holds  but $\mathfrak 
b=\omega_1$. Let $B=\{f_\alpha: \alpha<\omega_1\}$ be an unbounded subset of $\omega^\omega$. Since the filter of 
neighborhoods of $O$ in $V(\omega)$ is an $A$-filter, it follows that there exist $g_n\in \omega^\omega$, $n\in 
\omega$, such that every $F_{f_\alpha}$ contains $F_{g_n}$ for some $n\in \omega$; equivalently, for 
every $\alpha<\omega_1$, there exists an $n\in \omega$ such that $f_\alpha\leq^* g_n$. Let 
$g\in\omega^\omega$ be any function satisfying the condition $g_n\leq^* g$ for all $n\in\omega$. Then $f 
\leq^* g$ for all $f\in B$ and hence $B$ is bounded. This contradiction completes the proof. 
\end{proof}

It turns out that the existence of an $A$-filter on $\omega$ having no countable base is not only sufficient but 
also necessary for the existence of a nonmetrizable $\mathbb R^{\omega_1}$-factorizable space. 

\begin{theorem}
\label{Theorem_4.2}
The existence of an $A$-filter on $\omega$ which has no countable base is equivalent to the existence 
of a nonmetrizable (or, equivalently, non-second-countable) $\mathbb R^{\omega_1}$-factorizable space. 
\end{theorem}

\begin{proof}
In view of Corollary~\ref{Corollary_4.1}, it suffices to show that the existence of 
a non-second-countable $\mathbb R^{\omega_1}$-factorizable space implies that of an 
$A$-filter without a countable base on a countable set.  

Suppose that every $A$-filter on a countable set has a countable base. 
Let $X$ be an $\mathbb R^{\omega_1}$-factorizable space. First, note that $X$ is first-countable. Indeed, 
let $x\in X$ be a nonisolated point in $X$, and let $Y$ be a countable dense subspace of $X$ containing $x$. By 
Theorem~\ref{Theorem_3.3}\,(1) and Proposition~\ref{Proposition_4.1} 
the filter $\mathscr F_Y(x)$ of neighborhoods of $x$ in $Y$ is 
an $A$-filter and hence has a countable base $\mathscr B$. Since $X$ is regular and $Y$ is 
dense in $X$, it follows that $\{\overline B: B\in \mathscr B\}$ (the closures are in $X$) 
is a countable base of neighborhoods of $x$ in~$X$.

Thus, $X$ is first-countable. It obviously follows that the weight of any subspace $Z$ of $X$ of cardinality 
$\omega_1$ is at most $\omega_1$. Therefore, by Theorems~\ref{Theorem_3.2}\,(2) and~\ref{Theorem_3.3}\,(1) 
any such subspace $Z$ 
is second-countable. According to \cite[Theorem~1.1]{T_1978} (see also \cite{H-J}), $X$ is second-countable. 
\end{proof}

Note that if an $A$-filter $\mathscr F$ on $\omega$ has no countable base, then it cannot have a base of 
cardinality $\omega_1$ either, because the space $\omega_{\mathscr F}$ is $\mathbb
R^{\omega_1}$-factorizable (by Corollary~\ref{Corollary_4.1}) and  hence $w(\omega_{\mathscr F})>\omega_1$ 
(by Theorem~\ref{Theorem_3.2}\,(2)). However, if there exists an $A$-filter on $\omega$ which has no countable base, 
then there exists an $A$-filter on $\omega$ which has a base of cardinality exactly $\omega_2$ (and no base of 
smaller cardinality). Before proving this fact, we introduce some convenient notation. 

By analogy with filters, given families $\mathscr A,\mathscr B\subset [\omega]^\omega$, we say that $\mathscr B$ 
is a \emph{base} for $\mathscr A$ if every $A\in\mathscr A$ contains some $B\in\mathscr B$. 
Thus, a filter $\mathscr F$ is an $A$-filter if and only if every family $\mathscr A\subset\mathscr F$ 
with $|\mathscr A|\leq\omega_1$ has a countable base $\mathscr B\subset \mathscr F$. 

Given a family $\mathscr A\subset [\omega]^\omega$, we set 
\begin{align*}
I(\mathscr A) &=\{\textstyle\bigcap\gamma: \gamma\subset \mathscr A,\ |\gamma|<\omega\},\\ 
F(\mathscr A) &=\{L\subset\omega: \text{there exists an $M\in I(\mathscr A)$ such that $M\subset L$}\},\\ 
w(\mathscr A)&=\min\{|\mathscr B|: \mathscr B\subset F(\mathscr A),\ \mathscr B \text{ is a base for } 
F(\mathscr A)\}. 
\end{align*}
Note that $F(\mathscr A)$ is the smallest filter containing $\mathscr A$ and $\mathscr A$ is its subbase. Note 
also that if $\mathscr A$ is a filter, then $F(\mathscr A)=\mathscr A$. In this case, it is natural to refer to 
$w(\mathscr A)$ as the \emph{weight} of the filter~$\mathscr A$. Clearly, any filter of countable weight is an 
$A$-filter.

\begin{proposition}
\label{Proposition_4.5}
Any  $A$-filter of uncountable weight on $\omega$ contains an $A$-filter 
of weight $\omega_2$. 
\end{proposition}

\begin{proof}
Let $\mathscr F$ be an $A$-filter of uncountable weight on $\omega$. 
Then $w(\mathscr F)>\omega_1$ by Corollary~\ref{Corollary_4.1} and 
Theorem~\ref{Theorem_3.2}. 

\begin{lemma}
For any family $\mathscr A\subset \mathscr F$ with $|\mathscr A|\leq\omega_2$, there exists an 
$E(\mathscr A)\subset \mathscr F$ such that $\mathscr A\subset E(\mathscr A)$ and 
$w(E(\mathscr A))=| E(\mathscr A)|=\omega_2$. 
\end{lemma}
 
\begin{proof} 
If $w(\mathscr A)=\omega_2$, then we set $E(\mathscr A)=\mathscr A$. Suppose that $w(\mathscr 
A)<\omega_2$. We construct a family $\{M_\alpha:\alpha<\omega_2\}\subset \mathscr F$ completing $\mathscr A$ 
to the required family $E(\mathscr A)$ by induction on $\alpha<\omega_2$.

Take $\alpha<\omega_2$ and suppose that $\{M_\beta:\beta<\alpha\}\subset \mathscr F$ is already constructed. 
Let $\mathscr A_\alpha=\mathscr A\cup\{M_\beta:\beta<\alpha\}$. Since 
$w(\mathscr A_\alpha)\leq w(\mathscr A)+|\alpha|\le\omega_1$, it follows 
that $F(\mathscr A_\alpha)\neq \mathscr 
F$. Choose any $M_\alpha\in \mathscr F\setminus F(\mathscr A_\alpha)$. 

Having defined $M_\alpha$ for all $\alpha<\omega_2$, we set $\mathscr E=\mathscr 
A\cup\{M_\beta:\beta<\omega_2\}$. Note that $|\mathscr E|=\omega_2$. Let us show that $w(\mathscr 
E)=\omega_2$. Assume the contrary. Then $w(\mathscr E)<\omega_2$ and hence $F(\mathscr E)$ has a base $\mathscr 
B\subset F(\mathscr E)$ with $|\mathscr B|<\omega_2$. We have $\mathscr B\subset F(\mathscr A_\alpha)$ for 
some $\alpha<\omega_2$. Therefore, 
$$
F(\mathscr E)=F(\mathscr B)\subset F(\mathscr A_\alpha) \not\ni M_\alpha 
\in \mathscr E.
$$ 
This contradiction shows that $w(\mathscr E)=\omega_2$. We set $E(\mathscr A)=\mathscr E$. 
\end{proof}

\begin{lemma}
For any family $\mathscr A\subset \mathscr F$ with $|\mathscr A|\leq\omega_2$, there exists a $B(\mathscr 
A)\subset \mathscr F$ such that $\mathscr A\subset B(\mathscr A)$, $|B(\mathscr A)|\leq \omega_2$, and each 
$\mathscr C\subset \mathscr A$ with $|\mathscr C|\leq\omega_1$ has a countable base $\mathscr B\subset B(\mathscr 
A)$. 
\end{lemma}

\begin{proof}
Let us index (possibly with repetitions) the elements of $\mathscr A$ by ordinals: $\mathscr 
A=\{A_\alpha:\alpha<\omega_2\}$. For each $\alpha<\omega_2$, we set 
$\mathscr A_\alpha=\{A_\beta:\beta<\alpha\}$. 
Note that $w(\mathscr A_\alpha)\le\omega_1$. Since 
$\mathscr F$ is an $A$-filter, it follows that 
$\mathscr A_\alpha$ has a countable base $\mathscr B_\alpha\subset \mathscr F$ and that
the family $B(\mathscr A)=\mathscr A\cup
\bigcup\{\mathscr B_\alpha:\alpha<\omega_2\}$ is as 
required. 
\end{proof}

We proceed to prove the proposition. For each $\alpha<\omega_2$, we define 
a family $\mathscr A_\alpha\subset\mathscr F$ by induction as follows. We set $\mathscr A_0=\{\omega\}$. 
Assuming that $0<\alpha<\omega_2$ and $\mathscr A_\beta$ are already defined for all $\beta<\alpha$,  
we set 
$$
\mathscr A_\alpha^\ast=\bigcup_{\beta<\alpha}\mathscr A_\beta,\quad 
\mathscr E_\alpha=E(\mathscr A_\alpha^\ast),\quad \text{and}\quad 
\mathscr A_\alpha=B(I(\mathscr E_\alpha)).
$$

Having defined all $\mathscr A_\alpha$, we set 
$\mathscr A=\bigcup_{\alpha<\omega_2}\mathscr A_\alpha$ and $\mathscr 
G=F(\mathscr A)$. By construction, $|\mathscr A|=\omega_2$ and, therefore, $w(\mathscr G)\leq\omega_2$.

Let us show that $w(\mathscr G)=\omega_2$. Suppose that, on the contrary, $w(\mathscr G)<\omega_2$. 
Then $\mathscr G$ has a base $\mathscr B\subset \mathscr G$ with $|\mathscr B|\leq\omega_1$. Note that 
$\mathscr G=\bigcup_{\alpha<\omega_2}F(\mathscr A_\alpha)$. Therefore, $\mathscr B\subset 
F(\mathscr A_\alpha)\subset F(\mathscr E_{\alpha+1})\subset \mathscr G$ for some $\alpha<\omega_2$. 
Clearly, $\mathscr B$ is a base for $F(\mathscr E_{\alpha+1})$, which 
contradicts $w(\mathscr E_{\alpha+1})=\omega_2$.

It remains to show that $\mathscr G$ is an $A$-filter. Consider any $\mathscr M\subset \mathscr G$ with 
$|\mathscr M|\leq\omega_1$. From the same considerations as above, $\mathscr M\subset F(\mathscr E_\alpha)$ for some 
$\alpha<\omega_2$, and since $I(\mathscr E_\alpha)$ is a base for $F(\mathscr E_\alpha)$ and 
$|\mathscr M|\le \omega_1$, 
it follows that $\mathscr M$ has a base $\mathscr L\subset I(\mathscr E_\alpha)$ with 
$|\mathscr L|\leq\omega_1$. In turn, by the definition of $B(I(\mathscr E_\alpha))$, 
$\mathscr L$ has a countable base $\mathscr B\subset B(I(\mathscr E_\alpha))=\mathscr A_\alpha\subset \mathscr G$, 
as required. 
\end{proof}

\begin{remark}
\label{Remark_4.1}
There exist filters on a countable set which are not $A$-filters.
Indeed, it follows from Proposition~\ref{Proposition_4.5} that any filter of weight $\omega_1$ on a 
countable set is not an $A$-filter. To construct such a filter, it suffices to take an $\omega_1$-sequence of 
functions $f_\alpha\in \omega^\omega$, $\alpha<\omega_1$, such that 
$f_\alpha(n)<f_\beta(n)$  for all but finitely many $n\in \omega$ whenever $\alpha<\beta<\omega_1$ (it 
exists because any countable set of functions is bounded). The filter on $\omega\times \omega$ generated by the 
epigraphs $B_\alpha=\{(m,n)\in \omega\times \omega:  n\ge f_\alpha(m)\}$, $\alpha<\omega_1$, is as required.  
\end{remark}

\section{$\mathbb R^{\omega_1}$-Factorizable Groups}

We remind that a topological group $G$ is said to be \emph{$\mathbb R$-factorizable} if any continuous function  
$f\colon G\to \mathbb R$ factors through a continuous homomorphism to a second-countable group, 
that is, given any continuous function $f\colon G\to \mathbb R$, there exists 
a second-countable topological group $H$, a continuous homomorphism $h\colon G \to H$, and a continuous function 
$g\colon H\to \mathbb R$ for which $f = g \circ h$. Basic information about $\mathbb R$-factorizable groups can 
be found in \cite{A-T}. In particular, according to Proposition~8.1.3 of \cite{A-T}, any $\mathbb R$-factorizable 
group $G$ is \emph{$\omega$-narrow}, that is, given any neighborhood $U$ of the identity element of $G$, there 
exists a countable set $A\subset G$ for which $U\cdot A\subset G$. 

\begin{theorem}
\label{Theorem_5.1}
A topological group $G$ is $\mathbb R^{\omega_1}$-factorizable if and only if any continuous function  
$f\colon G\to \mathbb R^{\omega_1}$ factors through a homomorphism to a second-countable group, that is, given 
any continuous function $f\colon G\to \mathbb R^{\omega_1}$, there exists a second-countable topological group 
$H$, a continuous homomorphism $h\colon G \to H$, and a continuous function $g\colon H\to \mathbb R^{\omega_1}$ 
for which $f = g \circ h$. 
\end{theorem} 

\begin{proof}
The `if' part is obvious. Let us prove the `only if' part. Suppose that a topological group $G$ is $\mathbb 
R^{\omega_1}$-factorizable. 
From Theorem~\ref{Theorem_3.2}\,(1) it follows that $G$ is Lindel\"of and hence $\mathbb R$-factorizable \cite[Assertion~1.1]{T_1991}.
Let $f\colon G\to \mathbb R^{\omega_1}$ be a continuous function. By definition there exists a second-countable space 
$M$ and continuous functions $\varphi\colon G\to M$ and $\psi\colon M\to \mathbb R^{\omega_1}$ for which 
$f=\psi\circ \varphi$. Since $G$ is $\mathbb R$-factorizable, it follows from Lemma~8.1.2 of \cite{A-T} that there 
exists a second-countable group $H$, a continuous homomorphism $h\colon G\to H$, and a continuous function 
$\theta\colon H\to M$ for which $\varphi= \theta\circ h$. For $g=\psi\circ \theta$, we have $f= g\circ h$.
\end{proof}

\begin{remark}
Since any Tychonoff space $X$ with $w(X)\le \omega_1$ embeds in $\mathbb R^{\omega_1}$, it follows that a group 
$G$ is $\mathbb R^{\omega_1}$-factorizable if and only if any continuous function from $G$ to a Tychonoff space of 
weight at most $\omega_1$ factors through a homomorphism to a second-countable group. 
\end{remark}

In what follows, by $B(\omega)$ we denote the Boolean group generated by $\omega$,  
that is, the set $[\omega]^{<\omega}$ of all finite subsets of $\omega$ under the operation of symmetric 
difference (which we denote by +). The zero $\mathbf 0$ of $B(\omega)$ is the empty set $\varnothing \in 
[\omega]^{<\omega}$, and the inverse of each $x\in B(\omega)$ is $x$ itself. Given a free filter $\mathscr F$ 
on $\omega$, by $B(\mathscr F)$ we denote $B(\omega)$ with the group topology generated by the neighborhood base 
$\{\langle F\rangle : F\in \mathscr F\}$ at zero (here $\langle F\rangle $ denotes the subgroup generated by 
$F$; 
note that $\langle F\rangle =\{ g\in B(\omega): g\subset F\}$).

\begin{proposition}
\label{Proposition_5.1}
The group  $B(\mathscr F)$ is $\mathbb R^{\omega_1}$-factorizable if and only if $\mathscr F$ is a free 
$A$-filter. 
\end{proposition}

\begin{proof}
Let $\mathscr F_{\mathbf 0}$ denote the filter of neighborhoods of $\mathbf 0$.

Suppose that $\mathscr F$ is a free $A$-filter on $\omega$.
Since the topological group $B(\mathscr F)$ is countable and homogeneous,  
it suffices to check that $\mathscr F_{\mathbf 0}$ is an $A$-filter (by Proposition \ref{Proposition_4.2}).
Let $U_\alpha\in \mathscr F_0$, $\alpha<\omega_1$. 
For $\alpha<\omega_1$, there exists an $F_\alpha\in \mathscr F$ such that 
$\langle F_\alpha\rangle\subset U_\alpha$. Since $\mathscr F$ is an $A$-filter, 
there exist $G_n\in \mathscr F$, $n\in\omega$, such that each $F_\alpha$ contains $G_{n_\alpha}$ 
for some $n_\alpha<\omega$. We set $V_n=\langle G_n\rangle$, $n\in\omega$. 
Then $V_{n_\alpha}\subset U_\alpha$ for $\alpha<\omega_1$. 

Conversely, suppose that $B(\mathscr F)$ is $\mathbb R^{\omega_1}$-factorizable. 
By Proposition~\ref{Proposition_4.2} the 
filter $\mathscr F_{\mathbf 0}$  is an $A$-filter and, clearly, so is the filter $\{F\cap 
\omega: F\in \mathscr F_{\mathbf 0}\}$, which coincides with $\mathscr F$. 
\end{proof} 

\begin{theorem}
\label{Theorem_5.2}
The following conditions are equivalent:
\begin{enumerate}
\item[\rm(1)]
there exists an $A$-filter on $\omega$ which has no countable base;
\item[\rm(2)] 
there exists a nonmetrizable $\mathbb R^{\omega_1}$-factorizable space;
\item[\rm(3)] 
there exists a nonmetrizable $\mathbb R^{\omega_1}$-factorizable group.
\end{enumerate}
\end{theorem}
\begin{proof}
The implication (1)\,$\Rightarrow$\,(3) follows from Proposition~\ref{Proposition_5.1}, 
(3)\,$\Rightarrow$\,(2) is obvious, and 
(2)\,$\Rightarrow$\,(1) follows from Theorem~\ref{Theorem_4.2}.
\end{proof}

The following statements follow from Theorems~\ref{Theorem_4.1} and~\ref{Theorem_5.2} 
and Corollary~\ref{Corollary_4.2}.

\begin{corollary}
\begin{enumerate}
\item[\rm(1)]
If $\mathfrak b>\omega_1$, then there exists a nonmetrizable $\mathbb R^{\omega_1}$-factorizable group. 
\item[\rm(2)] 
If $2^\omega=\omega_1$, then any $\mathbb R^{\omega_1}$-factorizable group is second-countable. 
\end{enumerate}
\end{corollary}

In \cite{HPTZ} the notion of a $\tau$-fine topological group was introduced. An {\em $\omega_1$-fine} group is precisely a topological group in which the filter of neighborhoods of the identity element is an $A$-filter. 

\begin{theorem}
\label{Theorem_5.3}
A topological group is $\mathbb R^{\omega_1}$-factorizable if and only if it is $\omega$-narrow 
and $\omega_1$-fine.
\end{theorem}

\begin{proof}
The `only if' part follows from Proposition~\ref{Proposition_4.1}, Theorem~\ref{Theorem_3.2}, 
and the observation that any Lindel\"of topological group is $\omega$-narrow.

Let us prove the 'if' part. Let $G$ be an $\omega$-narrow $\omega_1$-fine topological group.
Since every $\omega$-narrow topological group embeds as a subgroup in a product of second-countable topological groups \cite[Theorem 3.4.23]{A-T}, we assume that $G$ is a subgroup of $\prod_{\iota\in I}G_\iota$, where each $G_\iota$ is a second-countable topological group.

According to Proposition~\ref{Proposition_3.2}\,(4), to prove that $G$ is $\mathbb R^{\omega_1}$-factorizable, 
it suffices to show that, given any $J\subset I$ with $|J|\le \omega_1$,
there exists an at most countable set $Q\subset I$ and a continuous map $\rho\colon \pi_Q(G)\to \pi_J(G)$
for which $\pi_J|_G = \rho\circ \pi_Q|_G$.

For ${\iota\in I}$, let $\{U_n(\iota):n\in\omega\}$ be a neighborhood base at the identity element of
$G_\iota$. We set $\mathscr U=\{G\cap \pi^{-1}_\iota(U_n(\iota)): \iota\in J, n\in\omega\}$.
Since $|\mathscr U|\leq\omega_1$ and $G$ is $\omega_1$-fine, it follows that 
there exists a countable family $\{V_n: n\in\omega\}$ of neighborhoods 
of the identity element of $G$ such that every  $U\in \mathscr U$ contains $V_n$  
for some $n\in\omega$. 
We can assume that $V_{n+1}\subset V_n$ for $n\in\omega$.
For each $n\in\omega$, there exists a finite set $Q_n\subset I$ and a neighborhood 
$W_n$  of the identity element in $\prod_{\iota\in Q_n}G_\iota$ such that 
$G\cap \pi_{Q_n}^{-1}(W_n)\subset V_n$. Let $Q=\bigcup_{n\in\omega}Q_n$.

We have $\ker \pi_Q|_{G}\subset \bigcap_{n\in\omega}V_n\subset \bigcap \mathscr U = \ker \pi_J|_{G}$. 
Given any neighborhood $U$ of the identity element in $\pi_J(G)$, there exists a neighborhood $W$ of the 
identity element in $\pi_Q(G)$ for which $\bigl(\pi_Q|_G\bigr)^{-1}(W)\subset \bigl(\pi_J|_G\bigr)^{-1}(U)$. 
Indeed, there exist $U_1,U_2,\dots, U_m\in \mathscr U$ such that $\bigcap_{i=1}^m U_i\subset \pi_J^{-1}(U)$ and 
an $n\in \omega$ such that $V_n\subset U_i$ for $i=1,2,\dots,m$; the set 
$W=\pi_Q(W_n\cap G)$ is as required. According to Proposition~1.5.10 of \cite{A-T}, there 
exists a continuous homomorphism 
$\rho\colon \pi_Q(G)\to \pi_J(G)$ for which $\pi_J|_G = \rho\circ \pi_Q|_G$.
\end{proof}

Under the assumption $\mathfrak{p}=\mathfrak{c}$, Malykhin \cite{Malykhin1979} 
(see also \cite[Theorem 4.5.22]{A-T}) constructed a maximal nonmetrizable countable Boolean group $G$.
In \cite{HPTZ}, it was noted that if $\omega_1<\mathfrak{p}$, then 
the group $G$ is $\omega_1$-fine.
Therefore, by Theorem~\ref{Theorem_5.3}, $G$ is $\mathbb R^{\omega_1}$-factorizable.
The proof that $G$ is $\omega_1$-fine was only briefly outlined in \cite{HPTZ}; 
below we give some details explaining why this is true.

The base of the topological group $G$ in question has the form  
$\bigcup_{\alpha<\mathfrak{c}}\tau_\alpha$, where each $G_\alpha=(G,\tau_\alpha)$ is metrizable and 
$\tau_\beta\subset \tau_\alpha$ for $\beta<\alpha<\mathfrak{c}$. 
Note that $G$ embeds in $\prod_{\alpha<\mathfrak{c}} G_\alpha$.

Let $J\subset \mathfrak{c}$, $|J|\leq\omega_1$. Since $\omega_1<\mathfrak{p}=\mathfrak{c}$ and 
the cardinal $\mathfrak{p}$ is regular \cite[Theorem 3.1]{vD}, it follows that 
there exists an $\alpha<\mathfrak{c}$ such that $\beta<\alpha$ for all $\beta\in J$. 
For $S=J\cup\{\alpha\}$, the space $\pi_S(G)$ is homeomorphic to $G_\alpha$ and hence second-countable. 
According to Proposition~\ref{Proposition_3.2}\,(3),  $G$ is $\mathbb R^{\omega_1}$-factorizable.

\section{Open Problems}

Although the existence of nonmetrizable $\mathbb R^{\omega_1}$-factorizable spaces is consistent with 
$\mathrm{ZFC}$, the following problem remains open. 

\begin{problem}
Is it true that any $\mathbb R^{\omega_1}$-factorizable space is cosmic, that is, has a countable network? 
\end{problem}

The next problem was already stated in the introduction.

\begin{problem}
What set-theoretic assumptions imply the existence of $A$-filters without a countable base on $\omega$?
\end{problem}

Results obtained above suggest the following two problems. 

\begin{problem}
Is it true that the existence of a strictly $\subset^*$-descending 
$\omega_2$-sequence of subsets of $\omega$ is equivalent to the existence of nonmetrizable 
$\mathbb R^{\omega_1}$-factorizable spaces? 
\end{problem} 

\begin{problem}
Does the negation $\lnot \mathrm{CH}$ of the continuum hypothesis  imply the existence of  nonmetrizable 
$\mathbb R$-factorizable spaces? 
\end{problem}

The $\mathbb R^{\omega_1}$-factorizability of topological groups is closely related to the $\mathbb 
R$-factorizability of products of groups. According to \cite{R-S_2025}, a product $G\times H$ of two 
$\mathbb R$-factorizable groups $G$ and $H$ is $\mathbb R$-factorizable in the sense of 
Definition~\ref{Definition_2.1} 
if and only if $G\times H$ is $\mathbb R$-factorizable as a topological group. In 
\cite{S_2022, S_2023, S_2025, R_2024} Lindel\"of (and hence $\mathbb 
R$-factorizable \cite[Assertion~1.1]{T_1991}) topological groups $G$ and $H$ for which $G\times H$ is not $\mathbb 
R$-factorizable were constructed. In \cite{S_2025} an example of Lindel\"of topological groups $G$ and 
$H$ such that $G\times H$ is not $\mathbb R$-factorizable, $G$ is second-countable, and $H$ has cellularity 
$2^\omega$ was given. This suggests the following problem.

\begin{problem}
Does there exist a separable metrizable topological group $G$ and a separable Lindel\"of topological group $H$ 
for which $G\times H$ is not $\mathbb R$-factorizable? 
\end{problem} 

In relation to products of topological groups, it is natural to introduce the following definition. 

\begin{definition}
We say that a topological group $G$ is \emph{productively $\mathbb R$-fac\-tor\-iz\-able} if the group $G\times H$ is 
$\mathbb R$-factorizable (as a topological group) for any $\mathbb R$-factorizable group~$H$. 

More generally, given a class $\mathscr C$ of topological groups, we say that a topological group $G$ is 
\emph{$\mathscr C$-productively $\mathbb R$-factorizable} if the group $G\times H$ is $\mathbb R$-factorizable 
(as a topological group) for any $\mathbb R$-factorizable group~$H\in\mathscr C$. 
\end{definition}

\begin{problem}
Describe all productively $\mathbb R$-factorizable topological \tolerance1000 groups. 
\end{problem}

\begin{problem}
Describe all $\mathscr C$-productively $\mathbb R$-factorizable topological groups for the classes 
$\mathscr C$ of topological groups with the 
following properties: 
\begin{itemize}
\item[(i)]
Lindel\"of;
\item[(ii)]
 \emph{ccc} Lindel\"of; 
\item[(iii)]
separable Lindel\"of;
\item[(iv)]
cosmic; 
\item[(v)]
countable; 
\item[(vi)]
metrizable countable. 
\end{itemize}
\end{problem}

In particular, the answer to the following question is unknown. 

\begin{problem}
Is it true that the product of any $\mathbb R$-factorizable group and any countable topological group is $\mathbb 
R$-factorizable? 
\end{problem}

Note that if the answer to this question is yes, then any $\mathbb R$-factorizable group is 
pseudo-$\aleph_1$-compact. Indeed, suppose that an $\mathbb R$-factorizable group $G$ is not 
pseudo-$\aleph_1$-compact. Let $\mathscr F$ be any filter on $\omega$. Then the group $B(\mathscr F)$ is 
countable, so that by assumption the group $G\times B(\mathscr F)$ is $\mathbb R$-factorizable (as a group). By 
Theorem~3.1 of \cite{R-S_2025} the product $G\times B(\mathscr F)$ is $\mathbb R$-factorizable in the sense of 
Definition~\ref{Definition_2.1}, and by Theorem~\ref{Theorem_3.1} \,$B(\mathscr F)$ is $\mathbb R^{\omega_1}$-factorizable. 
Therefore, according to Proposition~\ref{Proposition_5.1}, any filter $\mathscr F$ on a countable set must be an 
$A$-filter, which is not true by Remark~\ref{Remark_4.1}.  The problem of the existence of a 
non-pseudo-$\aleph_1$-compact $\mathbb R$-factorizable group was considered in detail in \cite{R-S_2025} (see also 
references cited therein).

On the other hand, a negative answer would solve another important problem of the theory of 
$\mathbb R$-factorizable groups, namely, that of preservation of the $\mathbb R$-factorizability of topological 
groups by continuous homomorphisms (see \cite[Problem~3.9]{T_1991a}, \cite[Problem~8.4.1]{A-T}, and 
\cite[Question~2]{R-S_2025}). Indeed, let $G$ and $H$ be topological groups such that $G$ is $\mathbb R$-factorizable, 
$H$ is countable, and $G\times H$ is not $\mathbb R$-factorizable (as a group). Let $H'$ be the topological group 
isomorphic to $H$ as an abstract group but endowed with the discrete topology. Since $G$ is $\mathbb 
R$-factorizable, it follows from Proposition~2.3 of \cite{R-S_2025} that the product $G\times H'$ is $\mathbb 
R$-factorizable in the sense of Definition~\ref{Definition_2.1}. By Theorem~3.1 of \cite{R-S_2025} it is also $\mathbb 
R$-factorizable as a topological group. It remains to note that $G\times H$ is the image of $G\times H'$ under a 
continuous homomorphism.

Similar problems can be stated for general topological spaces. 

\begin{definition}
Given a class $\mathscr C$ of topological spaces, we say that a topological space $X$ is 
\emph{$\mathscr C$-productively $\mathbb R$-factorizable} if the product $X\times Y$ is $\mathbb R$-factorizable 
for any space $Y\in \mathscr C$. 
\end{definition}

\begin{problem}
Describe all $\mathscr C$-productively $\mathbb R$-factorizable topological spaces for the 
class $\mathscr C$ of all topological spaces and for the classes $\mathscr C$ of topological spaces with the 
following properties: 
\begin{itemize}
\item[(i)]
pseudo-$\aleph_1$-compact;
\item[(ii)]
 \emph{ccc}; 
\item[(iii)]
separable; 
\item[(iv)]
separable metrizable; 
\item[(v)]
countable; 
\item[(vi)]
countable metrizable. 
\end{itemize}
The same problem is also interesting when $\mathscr C$ is the class of all topological groups, of topological 
groups with any of the properties \textup{(i)--(vi)}, and of $\mathbb R$-factorizable topological groups with any 
of the properties \textup{(i)--(vi)}. 
\end{problem}

The last four problems are related to $A$-filters.

\begin{definition}
Let $\mathscr F$ be a filter on $\omega$. We say that $\mathscr F$ is an $A^+$-filter (an $A^-$-filter) if, for 
any family $\mathscr G\subset \mathscr F$ with $|\mathscr G|\le\omega_1$, there exists an $F\in \mathscr F$ (an 
infinite $F\subset \omega$) such that $F\subset^* G$ for each $G\in \mathscr G$. 
\end{definition}

\begin{remark}
\label{Remark_6.1}
Recall that the cardinal $\mathfrak p$ is defined as the minimum cardinality of a family $\mathscr F\subset 
[\omega]^\omega$ such that the intersection of any finitely many elements of $\mathscr F$ is infinite but there 
exists no infinite set $A\subset \omega$ such that $A\subset^*F$ for all $F\in \mathscr F$ (see, e.g., 
\cite{vD}). It is easy to see that $\mathfrak p$ equals the minimum weight of a free filter on $\omega$ 
which is not an $A^-$-filter.
\end{remark}

Clearly, any $A^+$-filter is an $A$-filter and any $A$-filter on $\omega$ is an $A^-$-filter.

\begin{problem}
Does there exist an $A$-filter on $\omega$ which is not an $A^+$-filter?
\end{problem}

\begin{remark}
It follows from Remark~\ref{Remark_6.1} that if $\mathfrak p>\omega_1$, then any free filter of weight $\omega_1$ on 
$\omega$ is an $A^-$-filter but not an $A$-filter (in view of Remark~\ref{Remark_4.1}). A na\"\i ve example of an 
$A^-$-filter on $\omega$ which is not an $A$-filter can be constructed as follows. 

Let $N_0$ and $N_1$ denote the sets of even and odd nonnegative integers, respectively.  
Take the  Fr\' echet (cofinite) filter $\mathscr F_0$ on $N_0$ and any non-$A$-filter $\mathscr F_1$ 
on $N_1$ (see Remark~\ref{Remark_4.1}). The filter $\mathscr F=\{A\cup B:A\in \mathscr F_0,\ B\in \mathscr F_1\}$ on 
$\omega$ is an $A^-$-filter but not an $A$-filter. 
\end{remark}

However, we do not know the answer to the following question. 

\begin{problem}
Does there exist an ultrafilter on $\omega$ which is an $A^-$-filter but not an 
$A$-filter? 
\end{problem}

\begin{problem}
Is it true that under the assumption $\omega_1<\mathfrak p$ (or $\mathrm{MA} + \lnot \mathrm{CH}$) any 
ultrafilter on $\omega$ is an $A$-filter? 
\end{problem}

\begin{problem}
Find a $\mathrm{ZFC}$ example of an ultrafilter on $\omega$ which is not an $A$-filter. 
\end{problem}

\section*{Acknowledgments}
The authors are very grateful to the referee for carefully reading the manuscript, useful comments, and drawing 
their attention to the paper~\cite{HPTZ}.

\end{document}